\documentclass[12pt]{amsart}
\usepackage{amsmath,amssymb}

\DeclareMathOperator{\Ad}{Ad}

\DeclareMathOperator{\Aut}{Aut}

\newtheorem{question}{Question}

\newtheorem{theorem}{Theorem}[section]
\newtheorem{lemma}[theorem]{Lemma}
\newtheorem{proposition}[theorem]{Proposition}

\newtheorem{definition}[theorem]{Definition}
\newtheorem{remark}[theorem]{Remark}
\begin{document}

\newcommand\inv{^{-1}}
\newcommand\fh{\mathfrak h}
\newcommand\fg{\mathfrak g}
\newcommand\fa{\mathfrak a}
\title[Rigid and Almost Rigid Structures]{Nonexistence of invariant rigid structures 
 and invariant almost rigid structures}
\author[E.J.Benveniste and D.Fisher]{E. Jerome Benveniste and David Fisher}

\begin{abstract}
We prove that certain volume preserving actions of Lie groups and their lattices do not preserve
rigid geometric structures in the sense of Gromov.  The actions considered are the "exotic"
examples obtained by Katok and Lewis and the first author, by blowing up closed orbits in
the well known actions on homogeneous spaces.  The actions on homogeneous spaces all preserve
affine connections, whereas the action along the exceptional divisor preserves a projective structure.
The fact that these structures cannot in some way be "glued together" to give a rigid 
structure on the entire space is not obvious.  

We also define the notion of an almost rigid structure.  The paradigmatic example of a rigid structure
is a global framing and the paradigmatic example of an almost rigid 
structure is a framing that is degenerate along some exceptional divisor.  We show that the actions
discussed above do possess an invariant almost rigid structure.

Gromov has shown that a manifold with rigid geometric structure invariant under a topologically transitive
group action is homogeneous on an open dense set.  How generally this open dense set can be taken to be the 
entire manifold is an important question with many dynamical applications.  Our results indicate one way in 
which the geometric structure cannot degenerate off the open dense set.
\end{abstract}
\maketitle

\section{Introduction}

Let $G$ be a semisimple Lie group with all simple factors of real rank at least $2$ and
$\Gamma$ a lattice in $G$.  Zimmer and Gromov have conjectured that actions of $G$ or $\Gamma$ on
compact manifolds that preserve unimodular rigid geometric structures, for example 
an affine connnection and a volume form, should be "essentially classifiable" \cite{Gromov,Z1,Z2}.
However, there are "exotic" actions of these groups on compact manifolds where it had been unknown if an invariant 
rigid geometric structure existed. These "exotic" actions are on manifolds with non-linear fundamental
group and the actions admit non-trivial deformations, both of which cast doubt on the possibility of the
desired classification.
The purpose of this article is to show that these actions do not admit invariant rigid geometric structures.  
The simplest action we
consider is obtained by modifying the action of $SL_n(\mathbb Z)$ on ${\mathbb T}^n$.  In 
\cite{KL}, Katok and Lewis construct a new action from this action by "blowing up" the origin
in the algebro-geometric sense.  Throughout this article, in keeping with standard terminology from 
algebraic geometry, we will call the submanifold obtained
from the point by blowing up the exceptional divisor.  Katok and Lewis then further modify this action so that it is volume
preserving.  It is claimed in that paper that there is no invariant rigid geometric structure,
though no argument is given.  What is fairly clear is that the affine connection on ${\mathbb T}^n$
cannot be extended over the exceptional divisor.

In \cite{Benv}, the first author generalizes the Katok-Lewis construction to show that one can
construct other examples by blowing up closed orbits.  Let $G$ be a simple group, $H$ a semisimple
group and $\Lambda$ a lattice in $H$.  Then G acts on $H/{\Lambda}$.  Given appropriate algebraic 
conditions, there will be many closed $G$ orbits $X$ for this action, and one can "blow-up" the 
normal directions to $X$ along $X$.  (More generally, one could blow up along closed $G$ invariant 
submanifolds.)
For the specific examples considered in \cite{Benv}, it is shown that there is no invariant
affine connection.

Here we prove there is no rigid $A$-structure in the sense of Gromov on any of these examples.
The idea of the proof is very simple, and we outline it here in the case of the action on the 
blow-up of the torus.  This action agrees with the standard action on the torus on an open dense 
invariant set.  Assume there is an invariant rigid structure $\omega$.
We will use the invariance of the structure under $SL_n(\mathbb Z)$ to show that the structure 
is invariant under the local action of any (partially defined) affine vector field on the torus.
To prove this we use Iozzi's thesis.  To apply her result, we first lift the action to ${\mathbb R}^n$
(technically, ${\mathbb R}^n$ with the integers blown-up, but since we work on the open dense set that
is the complement of the integers, this is irrelevant).  Here the $SL_n(\mathbb Z)$ action extends
to a $SL_n(\mathbb R)$ action, and the content of Iozzi's theorem says exactly that $\omega$ is invariant
under $SL_n(\mathbb R)$.  (See theorem \ref{theorem:iozzi} for a precise statement.)  Since it 
is ${\mathbb Z}^n$ invariant
by construction, it will then be invariant under all of $SL_n(\mathbb R){\ltimes}{\mathbb R}^n$.  This action
does not descend to the torus, but it does descend locally, so we have a large set of vector fields that
preserve the structure locally.

The rest of the proof is conceptually simple, but somewhat computationally involved in practice.  We use 
one of the vector fields produced above to show that the structure cannot be rigid on at a point on the
exceptional divisor. 
We look at the infinitesimal isometry $\tau$ of the structure defined by a vector field $V$ and look at
what happens to $\tau$ as we tend towards the exceptional divisor.  By making an appropriate choice of $V$ we can
produce an infinitesimal isometry of order $k$ at a point on the exceptional divisor 
that is trivial up to order $k-1$.  This directly 
contradicts the definition of a rigid geometric structure. The vector field $V$ we construct would
vanish to order $k$ at the origin on the torus, were the exceptional divisor not present.  We show
by explicit computation that the "blowing-up" procedure reduces this order of vanishing by one.

In the other cases we discuss below, we follow essentially the same outline.  One can construct somewhat
more general examples, which we do not discuss here, since the main obstacle to a more general 
non-existence result is a generalization of Iozzi's theorem.

The paradigmatic example of a rigid geometric structure on $M$ is a global framing on some higher 
order frame bundle $P^k(M)$.  That such a framing is a rigid structure follows easily from the 
computation that a non-vanishing vector field on $\mathbb R$ is a rigid structure.  And
in fact, any rigid structure defines such a framing.  In the last section
of this article we define an {\em almost rigid $A$-structure} and show that the actions discussed
above do admit almost rigid $A$-structures.  The paradigmatic example of such a structure
is a global framing on $P^k(M)$ that degenerates to a subframing over some submanifold of $M$.  
The simplest such example is a vector field on $\mathbb R$ that vanishes at only finitely many points and only to 
finite order. 

Given a compact manifold $M$ and a rigid geometric structure $\omega$ on $M$.  Further assume that
the automorphism group of $(M,\omega)$ has a dense orbit.  One of the principal results of \cite{Gromov},
is that there is an open dense set $U$ in $M$ such that the local isometries of $\omega$ act 
transitively on $U$.  Since the set of local isometries are everywhere finite dimensional, this exhibits
a locally homogeneous structure on $U$ of the form $H/L$ where $H$ is a Lie group.  The proof of our results
involves explicitly identifying a large subgroup of $H$.  The proof of the results amount to showing that 
the homogeneous structure on $U$ cannot deteriorate too badly on $M{\backslash}U$.  Determing when $U=M$ is a
question with many potential dynamical applications.  It is, for instance, a key step in the celebrated work
of Benoist, Foulon, and Labourie \cite{BFL}. More interesting 
connections between dynamics and rigid structures
have recently been pointed out by Renato Feres \cite{F2}.
  
In section \ref{section:geometricstructures} we provide two examples of rigid geometric structures
on non-homogeneous manifolds.  For the first the manifold $N_1$ is ${\mathbb R}^n$ with the origin
blown up, and for the second the manifold $N_2$ is the compactification of $N_1$ obtained by viewing
the complement of the blow up as a subset of the real projective space $P^n$.  Both rigid structures
are invariant under natural actions of $SL_n(\mathbb R)$.  This paper grew out of an attempt
to use these structures to construct an $SL_n(\mathbb Z)$ invariant rigid structure on the
examples of Katok and Lewis.
Our proof can be interpreted as showing this to be impossible because the holonomy of the rigid
structure on the torus generates too many "extra" local isometries.  To illustrate this point,
we also exhibit a geometric structure on $N_2$ invariant under an action of 
$SL_n(\mathbb R){\ltimes}{\mathbb R}^n$ which preserves the exceptional divisor.  While the
tautological structure on $N_2$ for the $SL_n(\mathbb R)$ action is rigid, the tautological
structure for the $SL_n(\mathbb R){\ltimes}{\mathbb R}^n$ action is not.  This difference
is closely related to the fact that the second group has non-trivial unipotent radical, see
discussion at the end of section \ref{section:geometricstructures}.

We would like to thank Renato Feres for many helpful conversations in 
the early phases of this project.  Both authors would like to thank the Newton Institute
 for it's hospitality at the beginning of this project. The second author also thanks the NSF for funding provided
by grant number DMS-9902411.

\section{Actions}
\label{section:actions}

In \cite {KL}, Katok and Lewis describe certain actions of arithmetic
lattices in $SL(n, {\mathbb R})$. To fix notation, we will briefly
describe their construction. Let ${\mathbb T}^n$ denote the
$n$-dimensional torus, and let $p_+$ and $p_-$ denote the points
$(\frac{1}{2},0,...,0)$ and $(0,\frac{1}{2},0,...,0)$ on ${\mathbb T}^n$
(in terms of linear 
coordinates); then $p_+$ and $p_-$ are fixed points for the action of
the congruence subgroup $\Gamma_2 \subset SL(n, {\mathbb Z})$. Let $A \in
SL(n, {\mathbb Z})$ be such that $A p_+ = p_-$, 
let $l_+: \{v \in
{\mathbb R}^n | ||v|| < .5 \} \to {\mathbb T}^n $
be a linear coordinate system at $p_+$ and $l_- = A \circ l_+$ the
corresponding coordinates at $p_-$. Next let $S = ({\mathbb R}^n-\{0\}) / {\mathbb
R}^{\times}$ and for a vector $v \in {\mathbb R}^n$, denote its class in
$L$ by $[v]$.  Define 

$$
L = \{([v],u) \in S \times {\mathbb R}^n | u = yv {\rm \, for \, some \,} y
\in {\mathbb R} \}$$

\noindent
and let 
$$
q: L \to {\mathbb R}^n $$

\noindent
be defined by $q([v],u) = u$.
let $E \subset L$ be the set $\{([v],0) \}$ and $L^+$ (resp. $L^-$) be
the subsets 
$\{([v],u) \in S \times {\mathbb R}^n | u = yu {\rm \, for \, some \,} y > 0 \, ({\rm resp. \, y<0})  \}$. 

We will call $E$ the
exceptional divisor. Note that $SL(n,{\mathbb R})$ has a natural action on
$L$ and that $q$ is equivariant. 

Now form a new manifold $M$. Let $V = \{([v],u) \in L| ||u|| < 10^{-n} \}$,
and let $V^+$ (resp. $V^-$) be $V \cap L^+$ (resp. $V \cap L^-$). Let
$U = {\mathbb T}^n - \{p_1, p_2\}$ and form $M$ by gluing $U$ and $V$
together as follows: $l \in V^+$ (resp $\in V^-$) is glued to
$l_+(q(l))$ (resp $l_-(q(l))$). Define an action of $\Gamma $ on $M$ by
letting  $\Gamma$ act in the standard way on $U$ and on $E \subset
V$. It is immediate that this is a smooth action.

The action of $\Gamma$ on $M$ does not preserve a volume form. However,
Katok and Lewis describe a modified action that does. The necessary
modification is as follows: let $U$ and $V$ be as before, and let
$\rho: {\mathbb R}^n - \{0 \} \to {\mathbb R}^n - \{0 \}$ by $\rho (u) =
||u||^{- \alpha} u$, where $\alpha = \frac{n-1}{n}$. Form a manifold
$M^{\prime}$ by gluing $(l \in V^+$ (resp. $V^-$) to $l^+ (\rho
(q(l)))$ (resp. $l^- (\rho (q(l)))$). Define an action of $\Gamma$ on
$M^{\prime}$ by declaring it to be the standard actions on $U$ and $E
\subset V$. It is not difficult to check that this action is smooth
and preserves a smooth volume form.

We need expressions in local coordinates for the gluing maps from $V$
to $U$ that 
define $M$ and $M^\prime$. If $(X_1,...,X_{n-1}) \in {\mathbb R}^n - {0}$,
denote by $[X_1:...X_n]$ its class in $S$ (homogeneous
coordinates). Let $O \in E$; transforming by an element of $SL(n, {\mathbb
R})$ if necessary, we can assume that $O = [0:...:1]$. Then define
coordinates in a neighborhood of $O$ by setting
$$
x_i(P) = X_i/X_n, \; \; i = 1,...,n-1;
$$
$$
y (P)= q_n
$$
for $P =([X_1:...:X_n], (q_1,...,q_n))$; since $q_n X_i = q_i X_n$ for
$i = 1,...,n-1$, it follows that
\begin{equation}
\label{loc_glu}
q(P)= (x_1y,...,x_{n-1}y,y).
\end{equation}

\section{Geometric Structures}
\label{section:geometricstructures}

If $N$ is a manifold, we denote the $k$-th order frame bundle of $N$
by $F^k N$, and by $J^{s,k} N$ the bundle of $k$-jets at $0$ of
maps from ${\mathbb R}^s$ to $N$. If $N$ and $N^{\prime}$ are two
manifolds, and $f:N \to N^\prime$ is a map between them, then the
$k$-jet $j^k (f)$ induces a map $J^{s,k} N \to J^{s,k} N ^\prime$ for
all $s$.  We let $D^k (N)$ be the bundle whose
fiber $D^k _p$ at a 
point $p$ consists of the set of $k$-jets at $p$ of germs of
diffeomorphisms of $N$ fixing $p$. We abbreviate $D^k _0 ({\mathbb
R}^n)$ by $D^k$; this is a real algebraic group. We will need to
describe explicitly certain elements of $D^k$; for this purpose,
observe that we can represent each element uniquely, in terms of
standard coordinates $(\xi_1,...,\xi_n)$ on ${\mathbb R}^n$, in the form
$$
(P_1(\xi_1,...,\xi_n),...,P_n (\xi_1,..,\xi_n))
$$
where $P_1$, $P_2$,...,$P_n$ are polynomials of degree $\leq k$. We
denote the vector space of such polynomial maps of degree $\leq k$ by ${\mathcal
P}_{n,k}$. 

We will have occasion to make calculations with these objects in local
coordinates. Note that if $U \subset N$ is an open set, a choice of
coordinates on $c: U^\prime \subset {\mathbb R}^n \to U$ induces a
trivialization of $D^k U$:
$$
D^k U \to U^\prime \times D^k
$$
by 
$j^k _p (h) \to (c^{-1} (p), j^k _0 (h^\prime))$ where $h^\prime $ is
defined by: $h^\prime (v) = c^{-1}(h(c(c^{-1}(p) +v)))-c^{-1}(p)$.

The group $D^k$ has 
a natural 
action on $F^k (N)$, 
where $n$ is the dimension of $N$.
Suppose we are given an algebraic action of $D^k$ on a smooth algebraic
variety $Z$. Then following Gromov (\cite {Gromov}), we make the
following definition: 
\begin{definition}
\begin{enumerate}
\item An {\it A -structure} on $N$ (of order $k$, of type $Z$) is a smooth
map $\phi : F^k (N) \to Z$ equivariant for the $D^k$ actions.
\item With notation as above, the {\it $r$-th prolongation} of $\phi$,
denoted $\phi^r$, is the map $\phi^r : F^{k+r} N \to J ^{n,r}Z$
defined by $\phi^r = j^r(\phi) \circ \iota _k ^{r+k}$ where 
$\iota _k^{r+k} : F^{k+r} N \to J^{n,r} (F^k N)$ is the natural inclusion and 
$j^k(h) : J^{n,r} (F^k N) \to J^{n,r} Z$ is as before; this is an
A-structure of type $J^{n,r}Z$ and order $k+r$.
\end{enumerate}
\end{definition}
Equivalently, an A-structure of type $Z$ and order $k$ is a smooth section of
the associated bundle $F^k N \times _{D^k} Z$ over N.
Note that an A-structure on $N$ defines by restriction an A-structure
$\phi|_U$ on any open set $U \subset N$.
\begin{remark}
A-structures were introduced in \cite {Gromov}; a good introduction to
the subject, with many examples, can be found in \cite {Benoist}.  A 
comprehensive and accesible discussion the results of \cite{Gromov} concerning actions
of simple Lie groups can be found in \cite{F1}.
\end{remark}

\noindent
Note that if $N$ and $N^\prime$ are $n$-manifolds, and $h: N \to
N^{\prime}$ is a diffeomorphism, then $h$ induces a bundle map
$j^k(h): F^k N \to F^k N^\prime$. 

\begin{definition}
\begin{enumerate}
\item If $\phi: F^k N \to Z$, $\phi^\prime : F^k N ^\prime \to Z$  are
A-structures, a diffeomorphism $h: N \to N^\prime$ is an {\it isometry}
from $\phi$ to $\phi^\prime$ if $\phi ^\prime \circ j^k(h) = \phi$.
\item A {\it local isometry} of $\phi$ is a diffeomeorphism $h: U_1
\to U_2$, for open sets $U_1, U_2 \subset N$, which is an isometry
from $\phi|_{U_1}$ to $\phi|_{U_2}$.
\end{enumerate}
\end{definition}

\noindent
For $p \in M$ denote by $Is^{loc} _p (\phi)$ the pseudogroup of local
isometries of $\phi$ fixing $p$, and, for $l \geq k$, we denote by $Is
^l _p (\phi)$ the set of elements $j^l_p (h) \in D^l _p$ such that       
$j_p ^l (\phi \circ j^k _p(h)) = \phi^{l-k}$, where both sides are
considered as maps $F^{k+l} N \to J^{l-k}Z$. $Is _p ^l (\phi)$ is a
group, and there is a natural homomorphism $r_p ^{l;m}: Is _p ^l (\phi) \to 
Is_p^m (\phi)$ for $m < l$; in general, it is neither injective nor
surjective.

\begin{definition}
The structure $\phi$ is called $k$-rigid if for every point $p$, the
map $r_p ^{k+1;k}$ is injective.
\end{definition}

In order to provide some additional examples, as promised at the end of the
introduction, we recall the following lemma of Gromov.

\begin{lemma}
Let $V$ be an algebraic variety and $G$ a group acting algebraically on $V$.
For every $k$, there is a tautological $G$ invariant geometric structure 
of order $k$ on $V$, given by  $\omega:P^k(V){\rightarrow}P^k(V)/G$.  
This structure is rigid if and only if the action of 
$G$ on $P^k(V)$ is free and proper.
\end{lemma}

The conclusion in the first sentence is obvious.  The second sentence is
proven in section $0.4$, pages 69-70, of \cite{Gromov}.

{\bf Examples}

\noindent
\begin{enumerate}

\item The action of $G=SL_n(\mathbb R)$ on ${\mathbb R}^n$ is algebraic.
So is the action of $G$ on the manifold $N_1$ obtained by blowing up the
origin.  The reader can easily verify that the action of $G$ on $P^2(N_1)$ 
is free and proper.

\item We can compactify $N_1$ by $N_2$ by viewing the complement of the blow up as
a subset of the projective space $P^n$.  Another description of the same
action, which may make the rigid structure more visible to the naked eye, is
as follows.  $SL_{n+1}(\mathbb R)$ acts on $P^n$. Let $G$ be 
$SL_n(\mathbb R)<SL_{n+1}(\mathbb R)$ as block diagonal matrices with
blocks of size $n$ and $1$ and $1{\times}1$ block equal to $1$.  Then $G$
acts on $P^n$ fixing a point $p$.  We can obtain $N_2$ by blowing up the fixed point $p$.
The $G$ actions on both $P^n$ and $N_2$ are algebraic and again the reader
can verify that the action is free and proper on $P^2(N_2)$.

\item In the construction from $2$ above, there is an action
of a group $H$ where $H=SL_n(\mathbb R){\ltimes}{\mathbb R}^n$ and
$G=SL_n(\mathbb R)<SL_n(\mathbb R){\ltimes}{mathbb R}^n<SL_{n+1}(\mathbb R)$.
The $H$ action fixes the point $p$ and so also acts on $N_2$ 
algebraically.  However, over the exceptional divisor, the action
is never free on any frame bundle, since the subgroup ${\mathbb R}^n$
acts trivially to all orders at the exceptional divisor.

\end{enumerate}

The behavior in example $3$ above illustrates the fact that existence
of rigid structures is more complicated for algebraic groups which are
not semisimple, see the discussion in section $0.4.C.$ of \cite{Gromov}.  This
is related to our proof of the non-existence of rigid structures for
the Katok-Lewis examples, which depends heavily on the fact that any such
structure is locally invariant under the full affine group.

\section{An Application of Iozzi's Theorem}
\label{section:iozzi}

We shall use the following theorem, due to Iozzi:
\begin{theorem}
\label{theorem:iozzi}
{\bf (Iozzi)} Let $P$ be a principal bundle over a manifold $N$ with
structure group an algebraic group $H$ and let $Y$ be an algebraic $H$
space. Let $G$ be a simple Lie group with finite center that acts on
$P$ by principal bundle automorphisms, and suppose that for almost
every point $p \in N$ the stabilizer of $p$ in $G$ is noncompact. Let
$\Gamma \subset G$ be a lattice, and suppose that $s$ is a measurable
section of the associated bundle $P \times _H Y$ which is essentially
invariant under the action of $\Gamma$. Then it is essentially
invariant under the action of $G$.  
\end{theorem}

\noindent
If $p \in {\mathbb T}^n$, denote by $SA^{loc} _p$ the pseudogroup of local
diffeomorphisms fixing $p$ and preserving the standard connection and
volume form, and by $SA^k _p$ the group of $k$-jets at $p$ of elements
of $SA^{loc} _p$; a choice of special affine coordinates at $p$
induces an isomorphism from $SA^k _p$ to $SL(n, {\mathbb R})$.  

\begin{lemma}
Let $\Gamma \subset SL(n,{\mathbb Z})$ be a congruence subgroup, and let $U
\subset {\mathbb T}^n$ be a $\Gamma$ invariant open subset. Let $\phi: F^k
(U) \to Z$ be a $\Gamma$ invariant A-structure. Then at each point $p
\in U$, $SA^{loc} _p$ is locally contained in $Is_p ^{loc} (\phi)$ and
$SA^k _p \subset 
Is _p ^k (\phi)$. 
\end{lemma}

\noindent
By ``locally contained'' we mean the following: if $h: U_1 \to h(U_1)$
is a local diffeomorphism in $SA^{loc} _p$, there is an open set $U_1
^\prime$, $p \in U_1 ^\prime \subset U_1$, such that $h|_{U_1 ^\prime}
\in Is_p ^{loc} (\phi)$.

\noindent
{\bf Proof:} The second statement follows from the first. Now let
$\pi: {\mathbb R}^n \to {\mathbb T}^n$ be the natural 
projection and let ${\tilde U} = \pi^{-1} (U)$. Then $\phi$ defines a
$\Gamma$ invariant structure ${\tilde \phi}$ on ${\tilde U}$ which is
also invariant 
by the action of ${\mathbb Z}^n$ by translations. Then ${\tilde \phi}$ can
be thought of as a measurable section of $F^k {\mathbb R}^n \times _{D^k}
Z$ over ${\mathbb
R}^n$ by ignoring the null set ${\mathbb R}^n - {\tilde U}$. Since the action of $SL(n, {\mathbb R})$ on ${\mathbb R}^n$ has
non-compact stablilizers, the conditions of Iozzi's
Theorem are fulfilled, and ${\tilde \phi}$ is essentially invariant
under $SL(n, {\mathbb R})$. Since it is also invariant under ${\mathbb Z}^n$,
it is invariant under the special affine group $SL(n, {\mathbb R}) {\mathbb
R}^n$.

Now let $p \in U$ and let $h: U_1 \to h(U_1)$ be an element of $SA _p
^{loc}$ where $U_1, h(U_1) \subset U$ are neighborhoods of $p$, and
assume without loss of generality that $U_1$ is an evenly covered neighborhood
of $p$ for the covering map $\pi$; we want 
to show that there is a $U^\prime _1 \subset U_1$ such that
\begin{equation}
\label{loc_eq}
\phi
\circ j^k (h)|_{U_1 ^\prime} = \phi|_{U_1 ^\prime}
\end{equation}
$h$
lifts to a globally defined special affine transformation ${\tilde h}
: {\mathbb R}^n \to {\mathbb R}^n$, and by the preceding paragraph, we know
that ${\tilde \phi} \circ j^k ({\tilde h}) = {\tilde \phi}$ almost
everywhere on ${\mathbb R}^n$. Let ${\tilde U}_1$ be a connected component
of $\pi^{-1} (U)$, and let ${\tilde p}$ be the point in ${\tilde U}_1$
such that $\pi ({\tilde p}) = p$. Choose ${\tilde U}_1 ^\prime
\subset {\tilde U}_1$ be such that $h({\tilde U}_1 ^\prime) \subset
{\tilde U}$. Then since
$$
{\tilde \phi}|_{{\tilde U}_1 ^\prime}  \circ j^k
({\tilde h}) ={\tilde \phi}|_{h({\tilde U}_1 ^\prime)} \; \; a.e.
$$ 
and both sides are
continuous sections of $F^k {\mathbb R}^n \times _{D^k} Z$ over $h({\tilde
U}_1 ^\prime)$, the equality holds
everywhere on $h({\tilde U}^\prime _1)$; projecting down to
${\mathbb T}^n$ gives (\ref {loc_eq}).   

\section{The Katok-Lewis Examples}
\label{section:katoklewis}

\begin{proposition}
\label{proposition:kl}
There is no $\Gamma$-invariant rigid A-structure on the Katok-Lewis
example $M$.  
\end{proposition}  

\noindent
{\bf Proof:}  
Suppose that $\phi$ is a $\Gamma$-invariant A-structure of order
$k-1$; we shall show 
that it cannot be rigid. In fact, we shall show that if $O$ is a point
on the exceptional divisor $E$, there is a nontrivial element of the
kernel of the map $Is_O ^k (\phi) \to Is_O ^{k-1} (\phi)$.

In a neighborhood $W$ of $O$
on $E$, choose coordinates $(x_1,....,x_{n-1}, y)$, so that in terms
of these coordinates
\begin{enumerate}
\item $E
\cap W = \{y = 0 \}$, and $W^+ = _{def} W \cap L^+ = \{y > 0 \}$; 
\item  the gluing map from $W - E$ to $U$ is given by 
$$
(x_1,...,x_{n-1},y) \to (x_1y,...,x_{n-1}y,y);
$$
\end{enumerate}
where on $U$ we are using the coordinates $l^+$;
call this map $\mu$. All of our constructions will be expressed in
terms of these fixed coordinates.

We will define a family $\tau_p \in Is_p ^k (\phi)$ for $p \in W^+$,
where $\tau _p$ depends continuously on $p$. For $v \in {\mathbb R}^n$,
and $b \in {\mathbb R}$, define the map
$$
U(b,v): {\mathbb R}^n \to {\mathbb R}^n
$$
by 
$$
U(b,v)(v+(u_1,...,u_{n-1},y)) = v+(u_1+by, u_2,..,u_{n-1},y).
$$
By the lemma, we have
$$
\nu(b,p)= j_p ^k ( \mu^{-1} \circ U(b, \mu (p)) \circ \mu ) \in Is_p ^k (\phi)
$$
for any $b \in {\mathbb R}$ and $p \in W^+$. If we let $b:W \to {\mathbb R}$
be a continuous function, then $p \to \nu(b,p)$ defines a continuous
section of $D^k W$.

Recall that the choice of coordinates induces a trivialization
$$
\Xi : D^k W \to W \times D^k;
$$
write $\Xi (j) = (p, \sigma (j))$ for $j \in D_p ^k W$. Then $\sigma
(\nu(b,p) = j_0 ^k (H_{b,p})$ where $H_{b,p} : {\mathbb R}^n \to {\mathbb
R}^n$ is a polynomial map of degree $\leq k$ defined by
$$
H_{b,p}(\xi_1,...,\xi_{n-1},\eta) \equiv_k
$$
$$
\mu^{-1} ((x_1 + \xi_1)(y
+ \eta) + 
b \eta, (x_2 + \xi_2)(y + \eta),..., y +\eta)) - (x_1,...,x_{n-1},y)
$$  
$$
\equiv _k   ((\xi_1 + b \eta (y + \eta)^{-1}, \xi_2,..., 
\eta))
$$
$$
\equiv _k  ((\xi_1 + b \eta y^{-1}(1 - \eta/y + (\eta/y)^2 - ... +
(-1)^{k-1} (\eta/y)^{k-1}), \xi_2,..., \eta))
$$
(here, we have written $f \equiv _k g$ for $j^k _0 (f) = j^k _0 (g)$
where $f$ and $g$ are germs at 0 of maps ${\mathbb R}^n \to {\mathbb R}^n$.)
Now define a function $b : W \to {\mathbb R}$ by $b(x_1,...,x_{n-1}, y) =
y^k$; then $\lim _{p \to 0} H_{b(p),p)}$ exists in ${\mathcal P}_{n,k}$ and
equals $H$ where
$$
H(\xi_1,..,\xi_{n-1}, \eta) = (\xi_1 + (-1)^k\eta ^k, \xi_2, ..., \xi_{n-1},
\eta).
$$ 

Then if 
$$
\tau_p = \nu (b(p),p),
$$
$$
\lim_{p \to 0} \tau_p = \lim _{p \to 0} \Xi ^{-1} (p, j^k _0
(H_{b(p),p})) = \tau_0
$$
where $\tau_0 = \Xi ^{-1} (0, j_0 ^k (H))$.
Since $j_0 ^k (H)$ is in the kernel of the projection $D^k \to
D^{k-1}$, $\tau_0$ is in the kernel of the map $D^k _0 W \to D^{k-1}
_0 W$. On the other hand, since $\tau_p \in Is^k _p (\phi)$ for $p \in
W^+$, by continuity $\tau _0 \in Is _0 ^k (\phi)$; this violates
rigidity at $0$ and the proposition is proved.

Similarly, we have the following:
\begin{proposition}
There is no $\Gamma$-invariant rigid $A$-structure on $M^{\prime}$.
\end{proposition}

\noindent
{\bf Proof:} The argument is very similar to the one given
above. Let $U$, $V$, $E$ have the same meanings as before, except that

now they are to be viewed as subsets of $M^{\prime}$. Again, suppose
that $\phi$ is an invariant A-structure; suppose 
it is rigid of order $k$. By the
same argument as before, $Is_p ^k (\phi)$ contains $SA_p ^{loc}$ for all $p
 \in U$. Now choose a point $p$ in $E$, and choose coordinates
$(x_1,...,x_{n-1}, y)$ centered at $p$ on a neighborhood $W$ so that,
with respect to these
coordinates, the map $q: L \to {\mathbb R}^n$ has the form (\ref{loc_glu});
then in terms of $l^+$, the gluing map from $W^+ = _{def} W \cap L^+$
to $U$ is
\begin{equation}
\label{GLUING}
(x_1,...,x_{n-1},y) \to (x_1 ^2 y^2 + ... + x_{n-1} ^2 y^2 +
y^2)^{-\alpha/2} (x_1 y, x_2 y,..., x_{n-1} y, y).
\end{equation}
Now define new coordinates $(x_1 ^{\prime}, ..., x_{n-1} ^{\prime} , y
^{\prime})$ on $W$ by
$$
x^{\prime} _i = x_i, \; \; i=1,...,n-1;
$$
$$
y^\prime = y(1+x_1 ^2 + ... + x_{n-1} ^2)^{-\alpha/(2-2 \alpha)} =
y(1+x_1 ^2 +...+x_{n-1}^2)^{(1-n)/2}.
$$
Then in terms of these coordinates, the map (\ref {GLUING}) becomes
$$
(x_1 ^{\prime},..., x_{n-1} ^{\prime}, y^{\prime}) \to (x_1 ^{\prime}
y ^{\prime \; \delta},..., x_{n-1} ^{\prime} y ^{\prime \; \delta}, y
^{\prime \; \delta})
$$
where $\delta = 1- \alpha = \frac{1}{n}$. Call this map $\mu ^ {\prime}$. Let
$U(b,v)$ have the same meaning as in the
previous proposition; again, we conclude that 
$$
\nu^{\prime}(p, b) = j_p ^{k} (\mu ^{\prime \; -1} \circ U(b, \mu
^{\prime} (p)) 
\circ \mu ^{\prime})
$$
is in $Is ^k _p (\phi)$. Using the trivialization
$$
\Xi^\prime : D^k W \to W \times D^k 
$$
given by this choice of coordinates, we have that $\nu^\prime (p,b) =
\Xi ^{\prime -1} (p, j^k _0 (H^\prime _{b,p}))$ for a polynomial map
$H^\prime _{b,p} : {\mathbb R}^n \to {\mathbb R}^n$ defined by:
$$
H^\prime _{b,p} (\xi_1 ^{\prime}, ..., \xi_{n-1} ^{\prime} , \eta
^{\prime} ) \equiv_k
$$
$$
\mu^{\prime \; -1}((x_1 ^{\prime} +
\xi_1 ^{\prime})(y^{\prime} + \eta ^\prime)^\delta + b ((y ^{\prime} +
\eta ^{\prime})^\delta - y^{\prime \; \delta}),(x_2 ^{\prime} + \xi_2
^{\prime})(y^{\prime} + \eta ^\prime)^\delta,...
$$
$$
...,(x_{n-1} ^{\prime} +
\xi_{n-1} ^{\prime})(y^{\prime} + \eta ^\prime)^\delta, (y^{\prime} +
\eta ^\prime)^\delta ) - (x_1^\prime,...,x_{n-1} ^\prime, y ^\prime) 
$$
$$
\equiv _k (\xi_1 ^{\prime} + b(1- y^{\prime \; \delta}(y^{\prime} + \eta
^\prime)^{-\delta}) , x_2 ^\prime,...,\xi_{n-1} ^\prime, \eta
^{\prime})
$$
$$
\equiv _k (\xi_1 ^{\prime} + b(1- (1 + \eta
^\prime/y^{\prime})^{-\delta}) , x_2 ^\prime,...,\xi_{n-1} ^\prime, \eta
^{\prime})
$$

$$
\equiv _k (\xi_1 ^{\prime} + b(-r_1 \eta^\prime/y^\prime - r_2 (\eta ^\prime / y
^\prime)^2 - ... - r_k (\eta ^\prime / y^ \prime)^k
) , x_2 ^\prime,...,\xi_{n-1} ^\prime, \eta
^{\prime})$$

where $1 + r_1 X + r_2 X^2 + ...$ is the Taylor series of
$(1+X)^{-\delta}$. Then if we let $b(p) = y^{\prime \; k}$ as
before, we have
$$
\lim_{p \to 0} \nu^{\prime} (p, b(p)) = \tau^\prime _0
$$
where
$$
\tau ^\prime _0 (\xi_1^\prime, ..., \xi_{n-1}^\prime, \eta^\prime) = (\xi_1
^\prime + r_k \eta ^{\prime \; k}, \xi_2 ^\prime, ..., \xi_{n-1}
^\prime, \eta ^\prime).
$$
Again, continuity implies that $\tau ^\prime _0 \in Is _0 ^k (\phi)$,
which violates the $k$-rigidity of $\phi$ at 0.

\section{G actions}
\label{section:Gactions}
 
In this section we will show how to adapt the results of section \ref{section:katoklewis} 
to the setting of G actions. Here we will be considering examples where $G$ acts on $H/{\Lambda}$ 
and we blow up and glue along two closed $G$ orbits.

The proof follows the same outline as before.  First we use Iozzi's theorem to construct many
local vector fields that preserve the structure.  Second, we show that these vector fields
define infinitesimal isometries that degenerate on the exceptional divisor.

We first describe the adaptations necessary to use Theorem 4.1.  Here our model space is $H/{\Lambda}$ rather
than ${\mathbb T}^n$.  Given any point $p$ of $H$, we can define an $H$ action fixing $p$ by translating by $p{\inv}$,
acting by $Ad(h)$ and then translating back.  This defines a pseudogroup of local diffeomorphisms near any point $p$ of $H/{\Lambda}$, which
we denote by $H^{loc}_p$.  We can clearly restrict this action to any subgroup $G<H$, and
we denote the corresponding pseudogroup by $G^{loc}_p$. The version of Lemma 4.2. that we need is:

\begin{lemma}
\label{lemma:localaction}
Let $G$ and $H$ be simple real algebraic Lie groups, $G<H$.  Let $\Lambda$ be a lattice
in $H$, so $G$ acts on $H/{\Lambda}$.  Let $U{\subset}H/{\Lambda}$ be an open dense $G$ invariant set.
Let $\phi:F^k(U){\rightarrow}Z$ be a $G$ invariant $A$-structure.  Then at each point $p{\in}U$, $G^{loc}_p$
is locally contained in $Is_p^{loc}(\phi)$.  
\end{lemma}

\begin{proof}
We proceed as in the proof of Lemma 4.2.  Let $\pi:H{\rightarrow}H/{\Lambda}$ be the natural projection and let
${\tilde U}={\pi}{\inv}(U)$.  Now $\phi$ defines a geometric structure ${\tilde \phi}:{\tilde U}{\rightarrow}Z$
which is also right $\Lambda$ invariant.  We can think of $\phi$ as a measurable section of 
$F^k(H){\times}_{D^k}Z$.  Regretably here, the $\Lambda$ action does not have non-compact stabilizers, so we need
 an additional trick.  Since $\phi$ is $G$ invariant, we can think of $\tilde \phi$ as a $\Lambda$ invariant 
measurable section of the $Z$ bundle over $G{\backslash}H$ given by $G{\backslash}F^k(H){\times}_{D^k}Z$.
Since $G{\backslash}F^k(H){\rightarrow}G{\backslash}H$ is a $D^k$ bundle, and the stabilizers
of points for the $H$ action on $G{\backslash}H$ are all non-compact, we are in the setting where
Iozzi's theorem applies.  This implies that $\tilde \phi$ is invariant under the right $H$ action on 
$H$.  The rest of the proof of the lemma proceeds exactly as in Lemma 4.1.
\end{proof}

Now let $M$ be $H/{\Lambda}$ with a closed $G$ orbit blown-up as in \cite{Benv}.  Let $M'$ be the modification
of $M$ on which there is a $G$ invariant volume form as described in \cite{Benv}.  (The modification here
is analogous to the one described in section \ref{section:actions} for the Katok-Lewis examples and is 
described in detail in \cite{Benv}.) We have the following:

\begin{proposition}
\label{proposition:G}
There is no $G$ invariant rigid $A$-structure on either $M$ or $M'$.
\end{proposition}

In order to prove the proposition we must describe local coordinates on $M$ and $M'$
as in section \ref{section:actions}.  In fact, we will only prove the proposition for $M$ since the
modifications necessary for $M'$ are exactly as in section \ref{section:katoklewis}.

We outline some of the construction of $M$ from section 2 of \cite{Benv} in order to describe
the local coordinates in which we will carry out our computation.
Let $M_0=H/{\Lambda}$.  We can identify $TM_0$ with
$M_0{\times}{\fh}$.  The derivative action is given by $h(m,v)=(hm, Ad(h)v)$.  The tangent bundle
to the $G$ orbits is clearly isomorphic to $M_0{\times}{\fg}$.  Since $G$ is semisimple, there is a $G$
invariant complement to $\fg$ in $\fh$, which we will call $V$.  Now $M_0{\times}V$ is a $G$ invariant 
subbundle of $TM_0$.  Given a closed $G$ orbit $N$ in $M_0$, we can find a tubular neighborhood $U$ of $N$ and 
a neighborhood $U'$ in $V$ such that the natural map $exp_V:N{\times}U'{\rightarrow}U$ is a diffeomorphism.
We construct $M$ by blowing up along two closed $G$ orbits $N_1$ and
$N_2$ and gluing.  To do this we work in a tubular neighborhood
of each $N_i$.  Here we merely sketch the construction locally, the reader is referred to \cite{Benv} for details.

Let $S={\mathbb R}^m/{\mathbb R}^*$.  We take
a subset of ${\mathbb R}^n{\times}{\mathbb R}^m{\times}S$ of the form
$$L=\{(x,v,[u]){\in}{\mathbb R}^n{\times}{\mathbb R}^m{\times}S|v=cu {\rm \, for \, some \,} c{\neq}0 \}$$
and let $q:L^+{\rightarrow}{\mathbb R}^{n+m}$ be the canonical projection.  
Mimicking the construction in section \ref{section:actions}, we also introduce similar coordinates in a 
neighborhood of $N_2{\times}U_0$.  As before, we partition $L$ into three sets:
$$L^+=\{(x,v,[u]){\in}{\mathbb R}^n{\times}{\mathbb R}^m{\times}S|v=cu {\rm \, for \, some \,} c>0 \}$$
$$L^-=\{(x,v,[u]){\in}{\mathbb R}^n{\times}{\mathbb R}^m{\times}S|v=cu {\rm \, for \, some \,} c<0 \}$$ 
and $E=\{\vec{x}, 0, [u])\}$.   

Introduce local cooordinates $l^+$ on a neighborhood $N_1{\times}U_0$ of the form 
$( v_1, {\ldots}, v_m,x_1{\ldots}x_n)$ where the $x_i$ are coordinate along $U_0$ directions and
the $v_j$ are coordinates on $N_1$. Mimicking the construction in section \ref{section:actions},
we also introduce similar coordinates $l^-$ in a neighborhood of $N_2{\times}U_0$.  We form
$M$ by identifying taking $(M_0-{N_1{\cup}N_2}){\cup}L/{R}$ where $R$ is the relation
defined by identifying $q(L^+)$ with $l^+(N_1{\times}U_0-N_1)$ and similarly for 
$q(L^-)$ with $l^-(N_2{\times}U_0-N_2)$.  

We are being somewhat careless here.  The construction we are outlining does "blow-up" along
a subspace in a Euclidean space. Since $N_i$ are not contractible the maps
$l^+$ and $l^-$ above are not well-defined. To make sure that everything is canonically defined in a neighorhood
of $N_i$ in $M_0$ we should define $L$ intrinsically.  This is done in section 2 of \cite{Benv} and we 
refer the reader there for the necessary justifications.  Since we are only interested in a computation
done in the neighborhood of a point in $E$, the construction above is sufficient
for our purposes. 

\begin{proof}
We follow the proof of Proposition \ref{proposition:kl}.  The main difficulty is to replace the map
$U(b,v)$ used there by a local isometry on $H/{\Lambda}$ for which the computation does not become
intractible.  As in that proposition, assume there is a $G$ invariant $A$-structure $\phi$ of order $k-1$.
Once again, we will show that if $O$ is a point on the exceptional divisor $E$, there is a non-trivial
element in the kernel of the map $Is_O ^k (\phi) \to Is_O ^{k-1} (\phi)$.

In order to define the infinitesimal isometry and complete the computation here, we must chose our coordinates
and the vector field defining the infinitesimal isometry carefully.  We will use facts about the
structure of simple Lie groups and their Lie algebras, using \cite{Varad} as a reference.  We first
pick a Cartan subalgebra $\fa<\fh$, and then note that $\fh=\fa{\oplus}{\bigoplus}_{\alpha}X_{\alpha}$
where $X_{\alpha}$ is a root subspace for a root $\alpha$ and $\alpha$ runs over all roots.  We pick $\fa$
such that $\fa{\cup}\fg$ is a Cartan subalgebra of $\fg$. We chose
coordinates on $H$ that are obtained by exponentiating these coordinates so as to be able to use the 
Baker-Campbell-Hausdorff formula to describe multiplication.  We choose a vector field $V{\in}\fg$ such
that $V=X_{\alpha}$ for a simple root ${\alpha}$.  It follows that $[V,X_{\beta}]{\in}X_{\alpha+\beta}$ for
any root $\beta$
and that $[V,W]={\alpha(W)V}$ for $W{\in}\fa$.  We pick a vector $Y{\in}{\fh}$ such that $Y{\notin}\fg$
and such that $Y=X_\beta$ where $\beta$ is a simple root and $X_{\alpha+\beta}{\neq}{\emptyset}$.  Such a
$Y$ exists, since otherwise either all simple root spaces for $\fh$ are in $\fg$ forcing $\fg=\fh$ or $\fg$ commutes
with it's complement in $\fh$, in which case, $\fh$ cannot be simple. We further note that $X_{\alpha+\beta}$
is not in $\fg$.  Otherwise $X_{\beta}=[X_{-\alpha},X_{\beta+\alpha}]$ would also be in $\fg$.    
(For discussion of structure theory that makes all of this transparent, see section 4.3 of \cite{Varad}.)  

Letting $b$ be a real number, we observe that $\Ad(\exp(bV)(\exp(M))=\exp(M+[bV,M]+ b^2A(V,M))$ for any
$M$ in a small neighborhood of $0$ in $\fh$, where $A(X,M)$ is a sum of various higher order
brackets of $M$ with $V$ with coeffecients powers of $b$.  
This follows from the Baker-Campbell-Hausdorff formula, see e.g. section 2.15 of \cite{Varad}.  

We now choose our
coordinates $(v_1,{\ldots},v_m,x_1,{\ldots},x_{n-1},y)$ in a neighborhood $W$ of $0$ 
 subject to the following specifications.
The coordinates are given by exponentiating a basis for $\fh$ consistent with the decomposition $\fh=\fa{\oplus}{\bigoplus}_{\alpha}X_{\alpha}$.
We will let $y$ be the coordinate corresponding to $\exp(tY)$, $v_1$ the coordinate corresponding to $\exp(tV)$,
and $x_{1}$ the coordinate corresponding to $\exp(t[V,Y])$ which is nontrivial by the considerations
in the preceding paragraph.  We also assume that the $v_1$ coordinate corresponding to $\exp(V)$ is one.

Furthermore we choose the coordinates $(v_1, {\ldots}v_m,x_1,....,x_{n-1}, y)$, so that in terms
of these coordinates
\begin{enumerate}
\item $E
\cap W = \{x_1={\cdots}=x_m=y = 0 \}$, and $W^+ = _{def} W \cap L^+ = \{y > 0 \}$; 
\item  the gluing map from $W - E$ to $U$ is given by 
$$
(v_1,{\ldots},v_m,x_1,...,x_{n-1},y) \to (v_1,{\ldots},v_m, x_1y,...,x_{n-1}y,y);
$$
\end{enumerate}
where on $U$ we are using the coordinates $l^+$;
call this map $\mu$. All of our constructions will be expressed in
terms of these fixed coordinates.  

We will define a family $\tau_p{\in}Is_p^k(\phi)$ for $p{\in}W^+$,
where $\tau _p$ depends continuously on $p$. Let $w{\in}{\mathbb R}^{m+n}={\fh}$ and $V{\in}\fg$ be as above,
and $b{\in}{\mathbb R}$, define the map
$$
U(bV,w): {\fh} \to {\fh}
$$
by 
$$
U(bV,w)(w+(v_1,{\ldots}v_m,u_1,...,u_{n-1},y)) = 
$$
$$
w+Ad(\exp(bV))(v_1,{\ldots},v_m,u_1, u_2,..,u_{n-1},y).
$$
By lemma \ref{lemma:localaction}, we have
$$
\nu(bV,p)= j_p ^k ( \mu^{-1} \circ U(bV, \mu (p)) \circ \mu ) \in Is_p ^k (\phi)
$$
for any $b \in {\mathbb R}$ and $p \in W^+$. If we let $b:W \to {\mathbb R}$
be a continuous function, then $p \to \nu(b,p)$ defines a continuous
section of $D^k W$.

Recall that the choice of coordinates induces a trivialization
$$
\Xi : D^k W \to W \times D^k;
$$
write $\Xi (j) = (p, \sigma (j))$ for $j \in D_p ^k W$. Then $\sigma
(\nu(bV,p) = j_0 ^k (H_{bV,p})$ where $H_{bV,p} : {\mathbb R}^{m+n} \to {\mathbb
R}^{m+n}$ is a polynomial map of degree $\leq k$ defined by
$$
H_{bV,p}(\alpha_1,{\ldots}\alpha_m, \xi_1,...,\xi_{n-1},\eta) \equiv_k
$$

$$
\mu^{-1}(U(bV,\mu(p))(v_1+\alpha_1,{\ldots}, v_m+\alpha_m, (x_1 + \xi_1)(y
+ \eta),..., y +\eta))$$
$$- (v_1,{\ldots},v_m, x_1,...,x_{n-1},y)
\equiv_k
$$

$$
\mu^{-1}(v_1+\alpha_1+q_1(b),{\ldots}, v_m+\alpha_m+q_n(b), (x_1 + \xi_1)(y
+ \eta)+p_1(b),..., y +\eta))$$
$$ - (v_1,{\ldots},v_m, x_1,...,x_{n-1},y) \equiv_k
$$ 

$$
((\alpha_1+q_1(b),{\ldots},{\alpha_m+q_m(b)},{\xi_1 + p_1(b)}(y + \eta)^{-1}, \xi_2+p_2(b)(y + \eta)^{-1},..., 
\eta)
$$

\noindent
Here each $q_i$ is a polynomial in $b$ with zero constant term and linear term either zero or 
$\alpha_{j(i)}b$ for some $j(i){\in}\{1,{\ldots}m\}$.  Each $p_l$ for $l{\neq}1$ is a polynomial in 
$b$ with zero constant term and linear term either zero or $(x_{k(l)}{\eta}+\xi_{k(l)}(y+\eta))b$ for $k(l){\in}\{1,{\ldots},n-1\}$.  The remaining polynomial $p_1$
has zero constant term and has linear term ${\eta}b$.  This follows by direct computation using
the Baker-Campbell-Hausdorff formula as discussed above, the formula for $[V,M]$  deduced from the fact that 
$V$ is a root space for a simple root for $\fg$ and the fact that $Y$ is a rootspace for simple root space
for a simple root for $\fh$.  We note here that
any of the polynomials $q_i,p_j$ with trivial linear term is trivial, though we do not need this for
our computation. 

Let $p'$ be a point on the singular divisor, i.e, a point where $x_1={\cdots}=x_m=y=0$.  
Now define a function $b : W \to {\mathbb R}$ by $b(v_1,{\ldots}v_m, x_1,...,x_{n-1}, y) =
y^k$, make the substitution $(y+\eta){\inv}{\equiv}_k(1 - \eta/y + (\eta/y)^2 - ... +
(-1)^{k} (\eta/y)^{k})$ and compute as in the proof of Proposition \ref{proposition:kl}.
Then $\lim _{p \to p'} H_{b(p)V,p)}$ exists in ${\mathcal P}_{n,k}$ and
equals $H$ where
$$
H(\alpha_1,{\ldots},\alpha_m, \xi_1,..,\xi_{n-1}, \eta) = (\alpha_1,{\ldots},\alpha_m,\xi_1 + (-1)^k\eta ^k, \xi_2, ..., \xi_{n-1},
\eta).
$$ 

Then if 
$$
\tau_p = \nu (b(p),p),
$$
$$
\lim_{p \to p'} \tau_p = \lim _{p \to p'} \Xi ^{-1} (p, j^k _0
(H_{b(p),p})) = \tau_{p'}
$$
where $\tau_{p'} = \Xi ^{-1} (0, j_0 ^k (H))$.
Since $j_0 ^k (H)$ is in the kernel of the projection $D^k \to
D^{k-1}$, $\tau_{p'}$ is in the kernel of the map $D^k _0 W \to D^{k-1}
_0 W$. On the other hand, since $\tau_p \in Is^k _p (\phi)$ for $p \in
W^+$, by continuity $\tau _{p'} \in Is _{p'} ^k (\phi)$; this violates
rigidity at $p'$ and the proposition is proved.
\end{proof}
  
\section{Almost rigid structures}
\label{section:almostrigid}
In this section, we define a generalization of the rigid A-structures
and show that it is possible to construct such an object on the
Katok-Lewis examples. We continue to use the notation of section
\ref{section:actions}. 

\begin{definition}
An A-structure $\phi$ is called $(j,k)$-almost rigid (or just almost
rigid) if for every point $p$, $r^{k,k-1}_p$ is injective on the
subgroup $r^{k+j, k}(Is^{k+j}) \subset Is^{k}$.
\end{definition}

\noindent
Thus $k$-rigid structures are the $(0,k)$-almost rigid structures.

\noindent
{\bf Basic Example:} Let $V$ be an $n$-dimensional
manifold. Let $X_1, \ldots X_n$ be a collection of vector fields on
$M$. This defines an A-structure $\psi$ of type ${\mathbb R}^{n^2}$ on
$M$. If $X_1,{\ldots},X_n$ span the tangent space of $V$ at every point, then the
structure is rigid in the sense of Gromov. Suppose instead that there exists a 
point $p$ in $V$ and  
$X_1 \wedge \ldots \wedge X_n$ vanishes to order
$\leq j$ at $p$ in $V$.  
Then $\psi$ is  a $(j,1)$-almost rigid structure. Indeed,
let $p \in M$, and let $(x_1, \dots, x_n)$ be coordinates around
$p$. Suppose that in terms of these coordinates, $X_l =a_l ^m
\frac{\partial}{\partial x_j}$. Suppose that $f \in Is^{j+1}_p$.
We must show that $r_p ^{j+1,1}(f)$
is trivial. Let $(f^1, \ldots , f^n)$ be the coordinste functions of
$f$. Then $f \in Is^{j+1} _p$ implies that 
\begin{equation}
\label{1}
a_k ^l - a_k ^m \frac{\partial
f^l}{\partial x^m }
\end{equation}
 vanishes to order $j+1$ at $p$ for all $k$ and
$l$. Let $(b_k ^l)$ be the matrix so that $b_k ^m a_m ^l = \det(a_r
^s) \delta_k ^l$. Multiplying expression (\ref{1}) by $(b_k ^l)$, we see that
$\det(a_r ^s) (\delta_k ^l - \frac{\partial f^l}{\partial x^k})$
vanishes to order $j+1$. But since by assumption $\det(a_r ^s)$
vanishes to order $\leq j$, this implies that $(\partial f^l /
\partial x^k)(p) = \delta_k ^l$, so $r_p ^{j+1,1}(f)$ is the
identity, as required. 

If confused by the notation, the interested reader may find it enlightening to work out 
the basic example in the trivial case $n=1$.  Similar arguments can be given to show that
frames that degenerate to subframes are also almost rigid, provided the order of vanishing
of the form defining the frame is always finite.

Here we will show that the Katok-Lewis examples support an invariant
$(1,k)$-almost rigid structure. 
 
Observe that $GL(n, {\mathbb
R})$ acts on $J^2 ({\mathbb R}^n, {\mathbb R}^n, 0,0)$ by composition on the
right. 
Let $U \subset J^2
({\mathbb R}^n, {\mathbb R}^n,0,0)$ be the open subvariety of points where the
action is free, and let $Z$ be the quotient. Note that $D^k$ acts on
$V$ by composition on the left. 

\begin{definition}
The {\it canonical generalized connection} on $L$ is the equivariant
map $\Phi:F^2 M \to Z$ defined as follows: if $j \in F^2 _ z L$ let $f :
U \subset {\mathbb R}^n \to L$ be a map such that $j^2 _0 (f) = j$. Then
$\Phi (j)$ is the class in $V$ of $j^2 _0 (T_{f(p)} \circ q \circ
f) \in J^2 ({\mathbb R}^n, {\mathbb R}^n, 0,0)$.
\end{definition}

To justify this definition, we need to check that $j^2 _0 (T_{f(p)}
\circ q \circ f)$ is in the set $U$. This is clear if $p$ is not
in the exceptional divisor, for then this is the jet of a
diffeomorphism. Now suppose $p$ is in the exceptional divisor; we can
assume that $p = ([0:0: \cdots :0:1], (0, \cdots, 0))$. Then we can
choose coordinates 
$ \mu: {\mathbb R}^n \to L$ near $p$ as usual by $\mu (x_1, \cdots
,x_{n-1},y) = ([x_1:\cdots :x_{n-1}:1],(x_1y,\cdots x_{n-1}y,y))$ and
the map $q \circ \mu$ is $(x_1,\cdots, x_{n-1},y) \to (x_1y,\cdots
x_{n-1}y,y) $ and it is easy to check that the 
stabilizer of the 2-jet of this map for the left action of $GL(n, {\mathbb
R})$ is trivial.
\begin{proposition} 
The above construction defines a smooth (1,2)-almost rigid
A-structure, invariant for the $SL(n, {\mathbb
R})$ action on $L$.
\end{proposition}

\begin{remark}  \label {relation_to_connection} The structure $\Phi$
is the ``same as'' the standard 
connection on ${\mathbb R}^n$ on the complement of the exceptional divisor
in the following sense: The diffeomorphism $\pi : L - E \to {\mathbb R}^n
- 0$ allows one to pull back the standard connection to a connection
on $L- E$, which is the same as an equivariant map
$$
\Psi : F^2 L|_{L-E} \to D^2 (n,{\mathbb R})/GL(n, {\mathbb R}) \subset V;
$$
this map coincides with the restriction of $\Phi$ to $F^2|_{L-E}$.

In particular, if $h: U \subset L-E \to U^{\prime} \subset L-E$ is an
affine map, it is also an isometry of the $ Z$- structure.
\end{remark}

This construction allows one to put an almost rigid structure on the
manifold $M$, as follows. Recall that $M$ can be written as the union
of a subset $U \subset L$ and a subset $U' \subset {\mathbb T}^n$, where
the gluing map $h: U-(U \cap E) \to U'$ is affine. Let
$\Phi$ be the restriction of the canonical generalized connection on
$L$ to $U$. Remark \ref 
{relation_to_connection} shows that a connection defines a structure
$\Psi$ of
type $Z$ on $U'$, and that $h$ is an isometry of $Z$-structures. Thus
$\Phi$ and $\Psi$ paste together to form structure of type $Z$ on $M$.
Let us call this structure $\phi$; the construction shows that
it is invariant under the $\Gamma$ action on $M$. Now:

\begin{proposition}
The structure $\phi$ is (1,2)-rigid.
\end{proposition}
{\bf Proof:} We need to show that for every point $p$ in $M$,
$r^{2,1}_p: r^{3,2}(Is^3 _p (\phi))\to Is^1(\phi)$ is injective. This is clear
for $p$ not in $E$, since there $\phi$ defines a connection, which is
rigid in the usual sense. Now if $p \in E$, we may assume without
loss of generality that $p = ([0:0: \cdots :0:1], (0, \cdots, 0)$ and
that coordinates $(x_1, \ldots, x_{n-1},y)$ are chosen near $p$ so that $q$
is expressed in these coordinates as $(x_1, \ldots, x_{n-1},y) \to
(x_1y, \ldots ,x_{n-1}y, y)$. Now let $F = (F^1, \ldots, F^n) \in Is^3
_p (\phi)$; then
\begin{eqnarray*} \label {2}
F^1(x_1, \ldots x_{n-1}, y)F^n (x_1, \ldots, x_{n-1},y) &\equiv_3& x_1y \\
F^2(x_1, \ldots x_{n-1}, y)F^n (x_1, \ldots, x_{n-1},y) &\equiv_3& x_2y \\
&\vdots& \\
F^n(x_1, \ldots, x_{n-1},y) &\equiv_3& y.
\end{eqnarray*}
Now assume that $r^{3,1}(F)$ is trivial; then
\[ F^i (x_1, \ldots
x_{n-1} ,y) = x^i +P^i (x_1, \ldots ,x_{n-1}, y)
\]
 for $i=1, \dots ,n$,
where the $P^i$ vanish to order 1 at 0. Inserting these expressions
into (\ref 2), we see first, that $P^n$ vanishes to order 3 at 0;
then that $x_i P^n + y P^i + P^i P^n$ vanishes to order 3 at 0; but
since $P^n$ vanishes to order 3, this implies that $y P^i$ vanishes to
order 3, so that $P^i$ vanishes to order 2. This completes the proof.

\section{Gromov's theorem for almost rigid structures and other questions}
\label{section:speculation}

One of the principal results of \cite{Gromov} is:

\begin{theorem}\label{thm1f} Let $G$ be a simple Lie group and suppose that $G$ acts 
  analytically on a compact manifold $M$ preserving a
  volume and an analytic rigid geometric structure.  Further assume the action is ergodic. 
Then there exists a linear representation $\sigma :\pi_1(M) \to GL_n(\mathbb R)$ such
  that the Zariski closure of $\sigma (\pi_1(M))$ contains a group
  locally isomorphic to $G$.
\end{theorem}

At this juncture the question arises as to whether Gromov's theorem remains true 
for almost rigid structures or at least for the particular structure described above
on the Katok-Lewis example.  A careful examination of Gromov's proof shows that, except
at one step, he only uses the rigid structure to show that $M$ is locally homogeneous
on an open dense set.  The key step in proving that the manifold is locally homogeneous
is showing that infinitesimal isometries extend to local ones.  Gromov's proof of this
fact still works for almost rigid $A$-structures on $M$ that are rigid on an open dense subset
$U$ of $M$.  Since all the
examples above of almost rigid structures are rigid on an open dense set, this leads to
the natural:

\begin{question} Is a $(j,k)$-almost rigid structure 
on a manifold $M$ always $k$-rigid on an open dense set?  Is this true if we 
further assume there is a group action preserving the structure which has a dense orbit?
\end{question}

We believe the answer to this question is positive, even without the additional assumption.
The remaining obstacle to proving Gromov's theorem is surprisingly simple.  He proves that
a locally defined analytic Killing field $V$ of an analytic rigid $A$-structure on a simply
connected manifold $M$ has a unique analytic continuation.  This is an easy generalization of 
a similar theorem for $G$-structures of finite type by Amores \cite{A}. To complete the proof of Gromov's
theorem for almost rigid structures, one only needs to prove an analogue of Amores' theorem for
almost rigid $A$-structures.  However, here one needs to be more careful.

If we let $M$ be the modified Katok-Lewis example from above, we can show that no representation
exists for $N=(G{\times}M)/\Gamma$.  The proof is extremely indirect, and we merely sketch it here.
The fundamental group for $M$ is an $HNN$ extension of ${\mathbb Z}^n*{\mathbb Z}_2*{\mathbb Z}_2$
where the automorphism $T$ defining the $HNN$ extension switches the ${\mathbb Z}_2$'s.
The fundamental group of $N$ is ${\Gamma}{\ltimes}{\pi_1(M)}$ where $\Gamma<SL_n(\mathbb Z)$ is a subgroup
of finite index.  Note that $N$ is non-compact.
However, this is not an obstacle to a proof of Gromov's theorem and this exact difficulty is
dealt with in \cite{FZ}.  In that paper, it is also shown that if Gromov's representation $\sigma$
exists for $N$, it restricts to a non-trivial representation $\sigma$ of $\pi_1(M)$.  Using the explicit
construction of Gromov's representation, the fact that the
rigid structure on $M$ is the standard connection on ${\mathbb T}^n$ off the exceptional divisor $E$ 
and the analysis contained in \cite{FZ}, one can show that $\sigma(\pi_1(M)$ surjects onto $Z^n$.  
Furthermore, one can see that this surjection and $\sigma$ are both $\Gamma$-equivariant. The
$\Gamma$ action on $\pi_1(M)$ comes from the inclusion $\Gamma<\Aut(\pi_1(M))=SL_n(\mathbb Z)$  
If we now apply the main theorem
of \cite{FW}, we see that this implies that there is a $\Gamma$ equivariant map $\phi$ from $M$ to
${\mathbb T}^n$ inducing the surjection on fundamental groups. 
It is easy to see that this is impossible by showing that $\phi$ is the identity off the exceptional
divisor and then considering the possibilities for $\phi(E)$.

The argument of the preceeding paragraphs is clearly unsatisfactory as a proof that analytic
continuation of vector fields is not always possible for almost rigid $A$-structures.  We note here
that the above argument fails to produce a contradiction if one simply passes to an
appropriate double cover of $M$ before beginning the construction.  It would be interesting to 
determine a natural additional condition on an almost rigid structure that would allow one to prove
an analogue of Amores' result.

\end{document}